\documentclass[12pt,twoside,reqno]{amsart}
\usepackage{amsmath}
\usepackage{amsrefs}
\usepackage{graphicx}
\usepackage{amsfonts}
\usepackage{amssymb}
\usepackage{amsthm}
\usepackage[mathscr]{eucal}
\usepackage{verbatim}
\usepackage{epsfig}
\usepackage{fmtcount}
\usepackage{cite}
\usepackage{upgreek}

\def\BB{{\mathcal{B}}}

\def\LL{{\mathcal{L}}}

\def\XX{{\mathcal{X}}}

\def\UU{{\mathcal{U}}}
\def\R{{\mathbb R}}
\def\C{{\mathbb{C}}}

\def\d{\delta}

\def\e{\varepsilon}
\def\f{\varphi}
\def\s{\sigma}
\def\a{\alpha}
\def\g{\gamma}
\def\l{\lambda}
\def\L{\Lambda}
\def\k{\varkappa}
\def\O{\Omega}

\def\XXint#1#2#3{{\setbox0=\hbox{$#1{#2#3}{\int}$}
\vcenter{\hbox{$#2#3$}}\kern-.5\wd0}}

\newtheorem{theorem}{Theorem}[section]
\newtheorem{lemma}[theorem]{Lemma}
\newtheorem{corollary}[theorem]{Corollary}

\newtheorem*{definition}{Definition}
\numberwithin{equation}{section}
\newcommand{\supp}{\mathop{\rm supp}\nolimits}

\DeclareMathOperator{\dist}{dist}

\begin{document}
\title[Singular operators with antisymmetric kernels]
{Singular operators with antisymmetric kernels, related capacities, and Wolff potentials}
\author{D.~R.~Adams}
\address{Department of  Mathematics, University of Kentucky, Lexington, KY 40506}
\email{dave@ms.uky.edu}
\author{V.~Ya.~Eiderman}
\address{Department of  Mathematics, University of Wisconsin-Madison, Madison, WI 53706}
\email{eiderman@math.wisc.edu}
\begin{abstract}
We consider a generalization of the Riesz operator in $R^d$ and obtain estimates for its
norm and for related capacities via the modified Wolff potential. These estimates are
based on the certain version of $T1$ theorem for  Calder\'on--Zygmund operators in metric
spaces. We extend two versions of Calder\'on--Zygmund capacities in $R^d$ to metric
spaces and establish their equivalence (under certain conditions). As an application, we
extend the known relations between $s$-Riesz capacities, $0<s<d$, and the capacities in
Nonlinear Potential Theory, to the case $s=0$.
\end{abstract}
\maketitle

\section{Introduction}%Section 1

For $\e>0$, $s>0$, and for a finite Borel measure $\mu$ on $\R^d$, $d\ge1$, define the $\e$-truncated $s$-Riesz transform of $\mu$ by the equality
$$
R_{\mu,\e}^s(x)=\int_{|y-x|>\e}\frac{y-x}{|y-x|^{s+1}}\,d\mu(y),\quad x,y\in\R^d.
$$
For $0<s\le d$ the limit
$$
R_{\mu}^s(x)=\lim_{\e\to0+}R_{\mu,\e}^s(x)
$$
exists almost everywhere in $\R^d$ with respect to Lebesgue measure; this limit is said to be the $s$-Riesz transform of $\mu$ at $x$. Analogously, we define the $\e$-truncated $s$-Riesz operator
$$
\mathfrak{R}_{\mu,\e}^s f(x)=\int_{|y-x|>\e}\frac{y-x}{|y-x|^{s+1}}\,f(y)\,d\mu(y),\quad f\in L^2(\mu),\quad \e>0.
$$
Then for every $\e>0$, the operator $\mathfrak{R}_{\mu,\e}^s$ is
bounded on $L^2(\mu)$. We set
$$
\pmb|\mathfrak{R}_{\mu}^s\pmb|:=\sup_{\e>0}\|\mathfrak{R}_{\mu,\e}^s\|_{L^2(\mu)\to L^2(\mu)}.
$$

It is known (see \cite{P}, \cite{MPV}, and \cite{ENV}) that for $0<s<1$
\begin{equation}\label{E:Wolff_energy}
\int_Q |\mathfrak{R}_{\mu,\e}^s \chi_Q|^2 d\mu \approx \int_Q
W^{\mu\mid Q} d\mu,
\end{equation}
where $Q$ is a cube in $\R^d$ with sides parallel to the coordinate axes, and $W^{\mu}$ is the acclaimed Wolff potential from Nonlinear Potential Theory; see \cite{HW} and \cite{AH}:
\begin{equation} \notag
W^{\mu}(x) = \int_0^{\infty} \left( \frac{\mu(B(x,t))}{t^s}\right)^2 \frac{dt}{t}.
\end{equation}
The symbol $\approx$ in (\ref{E:Wolff_energy}) means that the ratio is bounded above and below by positive constants that depend only on $s$ and $d$.  In the language of \cite{AH}, $W^{\mu}$ corresponds to the case $\frac{1}{p-1} = 2$ or $p = 3/2$ and $s = N - \alpha \frac{3}{2}$; see page 110 of \cite{AH}.  The right side of (\ref{E:Wolff_energy}) is consequently called the Wolff energy. The upper bound in (\ref{E:Wolff_energy}) holds for all $s \in (0,d)$ \cite{ENV}; the lower bound appears in \cite{P} -- though only for $0<s<1$. The latter is not correct when $s=0,1,\dots,d$. The important question about the validity of the lower bound for non-integer $s\in(1,d)$ remains open.

Relation (\ref{E:Wolff_energy}) plays a key role in the uniform boundedness of $\mathfrak{R}_{\mu,\e}^s$: (\ref{E:Wolff_energy}) holding for all cubes $Q$ and the non-homogeneous $T1$ theorem by Nazarov-Treil-Volberg  \cite{NTV97}, \cite{NTV2002}, \cite{NTV} implies that $\pmb|\mathfrak{R}_{\mu}^s\pmb|^2\le C\supp_{x\in\supp\mu}W^{\mu}(x)$ (see \cite[Theorem~2.6]{ENV}). Also, exploiting the connection with Non-Linear Potential Theory, this estimate implies that
\begin{equation}\label{E:second}
\gamma_{s,+}(E) \geq c \cdot \dot{C}_{\frac{2}{3}(d-s), \frac{3}{2}} (E)
\end{equation}
for any compact set $E \subset \R^d$, $0 < s < d$  \cite{ENV}, and
\begin{equation}\label{E:third}
\gamma_{s,+}(E) \leq c' \cdot\dot{C}_{\frac{2}{3}(d-s), \frac{3}{2}} (E)
\end{equation}
whenever $0<s<1$ \cite{MPV}.  Here $\dot{C}_{\alpha, p}(\cdot)$ is the Riesz capacity of order $\alpha$ and degree $p$ from \cite{AH}:
$$
\dot{C}_{\a,p}(E)=\sup_{\mu\in M_+(E)}\bigg(\frac{\mu(E)}{\|I_\a\ast\mu\|_{p'}}\bigg)^p,\quad I_\a(x)=\frac{A_{d,\a}}{|x|^{d-\a}},\quad \frac1{p'}+\frac1p=1,
$$
where $1<p<\infty$, $0<\a p<d$, $\|\cdot\|_{p'}$ is the $L^{p'}$-norm with respect to the Lebesque
measure in $\R^d$, and $A_{d,\a}$ is the certain constant depending on $d$ and $\a$; furthermore,
$$
\g_{s,+}(E):=\sup\{\|\mu\|:\mu\in M_+(E),\ \|R_{\mu}^s(x)\|_\infty\le1\},
$$
where $M_+(E)$ is the class of positive Radon measures supported on $E$.
The study of these set functions has accelerated recently with the breakthrough results of X.~Tolsa and others. In particular, Tolsa proved that $\gamma_{1,+}(\cdot)$ is comparable with the analytic capacity when $d=2$.

Clearly, $(I_\a\ast\mu)(x)\approx\|\mu\|\cdot|x|^{\a-d}$ for a finite measure $\mu$ with compact support and for sufficiently big $|x|$. If $\a p=d$, then $p'=d/(d-\a)$, and we see that $I_\a\ast\mu\not\in  L^{p'}$. Because the case $\a p=d$ will be important, we consider the standard Bessel capacity instead of the Riesz one, defined in the similar way:
$$
C_{\a,p}(E)=\sup_{\mu\in M_+(E)}\bigg(\frac{\mu(E)}{\|G_\a\ast\mu\|_{p'}}\bigg)^p,\quad 1<p<\infty,\quad \frac1{p'}+\frac1p=1,
$$
where $G_\a$ is the Bessel kernel. We refer to \cite{AH}, p.~9--13, for definitions and properties of the Bessel kernel and Bessel potentials. It is important to note that $G_\a(x)\approx I_\a(x)$ as $|x|\to0$, $0<\a<d$, and $G_\a(x)=O(e^{-c|x|})$ as $|x|\to\infty$, $\a>0$, $0<c<1$. Thus, $\dot{C}_{\a,p}(E)\approx C_{\a,p}(E)$ for compact sets $E$ with diam$(E)\le1$ and for $1<p<\infty$, $0<\a p<d$.

Inequality (\ref{E:third}) does not hold, in general, for integer $s\in(0,d]$. Indeeed, a closed $s$-dimensional ball is an example of a compact set with positive $\gamma_{s,+}$-capacity and zero $\dot{C}_{\frac{2}{3}(d-s), \frac{3}{2}}$-capacity. For $s=0$,  (\ref{E:third}) also does not hold. Notice that $\gamma_{0,+}(E) \geq1$ for every set $E$ in $\R^d$, but $\dot{C}_{\frac{2}{3}d, \frac{3}{2}} (\overline{B}(0,r))\to0$ as $r\to0$; here $B(x,r):=\{y\in\R^d:\ |y-x|<r\}$. The validity of (\ref{E:third}) for non-integer $s\in(0,d)$ is an open question; essentially it is equivalent to the problem about the lower bound in \eqref{E:Wolff_energy}.

This note is inspired by the following question: is there a natural analog of the capacity $\gamma_{s,+}$ which is equivalent to $C_{\frac{2}{3}(d-s), \frac{3}{2}}$ when $s$ is an integer? The particular interest is the case $s=0$. To be more precise, we generalize the notion of $\gamma_{s,+}$ in the following way. Let $\f(t)$ be a continuous increasing function of $t\ge0$ with $\f(0)=0$. We define the $\e$-truncated $\f$-operator and the $\f$-transform by the equalities
$$\aligned
\mathfrak{R}_{\mu,\e}^\f f(x)&=\int_{|y-x|>\e}\frac{y-x}{|y-x|}\cdot\frac1{\f(|y-x|)}\,f(y)\,d\mu(y),\quad f\in L^2(\mu),\quad \e>0,\\
R_{\mu,\e}^\f(x)&=\mathfrak{R}_{\mu,\e}^\f {\bf1}(x),
\quad R_{\mu}^\f(x)=\lim_{\e\to0}R_{\mu,\e}^\f(x),\quad x,y\in\R^d.
\endaligned$$
We assume that the limit exists almost everywhere in $\R^d$ with respect to Lebesgue measure. As above we set
$$
\pmb|\mathfrak{R}_{\mu}^\f\pmb|:=\sup_{\e>0}\|\mathfrak{R}_{\mu,\e}^\f\|_{L^2(\mu)\to L^2(\mu)}.
$$
The $\f$-Wolff potential of a Borel measure $\mu$ is defined by the formula
$$
W_\varphi^\mu (x) = \int_0^\infty \left(\frac{\mu(B(x,t))}{\varphi(t)}\right)^2 \frac{d\varphi(t)}{\varphi(t)}.
$$
When $\f(t)=t^s$, we write $R_{\mu,\e}^s$ instead of $R_{\mu,\e}^\f$, etc.

Let $\Sigma_\f$ be the class of nonnegative Borel measures $\mu$ in $\R^d$ such that
\begin{equation}\label{f14}
\mu(B(x,r))\le \f(r) \quad\text{for all}\ x\in\R^d\text{ and }r>0.
\end{equation}
We introduce the capacity $\gamma_{\f,+}$ of a compact set $E$ in $\R^d$ in the similar way, namely
$$
\g_{\f,+}(E):=\sup\{\|\mu\|:\mu\in\Sigma_\f,\ \supp\mu\subset E,\ \|R_{\mu}^\f(x)\|_\infty\le1\}.
$$
The condition $\mu\in\Sigma_\f$ in this definition is superfluous for $\f(t)=t^s$. Namely, it is shown in \cite{MPV}, p.~217, that if $\|R_{\mu}^s(x)\|_\infty\le1$, $0<s<d$, then
$\mu(B(x,r))\le Cr^s,\ x\in \R^d,\ r>0$, for every measure $\mu\in M_+(E)$.  For $s=d-1$, this fact is also noted in \cite{V}, p.~46. We do not know if the condition $\mu\in\Sigma_\f$ can be droped for any $\f$.

Now we can formulate our question as follows. Given integer $s\in[0,d]$, is it possible to find $\f$ for which
\begin{equation}\label{f15}
\gamma_{\f,+}(E) \approx C_{\frac{2}{3}(d-s), \frac{3}{2}} (E)
\end{equation}
for all (or at least for sufficiently small) compact sets $E \subset \R^d$? With this end in view we extend some results of \cite{ENV} to the class of Calder\'on-Zygmund (CZ) operators on separable metric spaces. As an application, we obtain an extension of results in \cite{MPV} and \cite{ENV} to the class of concave and convex functions $\f$ satisfying the doubling condition. In particular, we derive the theorem on the comparison of the capacity $\gamma_{\f,+}$ and $\f$-Wolff potentials. These generalizations are, we believe, of independent interest. As a corollary, we give an affirmative answer the question posed above, for $s=0$, and indicate the corresponding function $\f$. Namely, we prove that
$\gamma_{\f,+}(E) \approx C_{\frac{2}{3}d, \frac{3}{2}} (E)$ for every compact set $E$ with diam$(E)\le1$, if $\f$ is a concave increasing function on the interval $(0,\infty)$ such that
$\varphi(t) = \varphi_0(t) =(\log \frac{1}{t})^{-1/2}$ as $0 < t\le  e^{-3/2}$,  and $\f(2t)\le2^s\f(t)$, $0<t<\infty$, with some $s\in(0,1)$.

We conjecture that for positive integers $s$ the answer is negative.

By $c,C$ we denote various positive constants.

\section{Main results}%Section 2

We start with introducing of the class of functions $\f$.
\begin{definition} By $\Phi$ we denote the class of functions $\f(t),\ t\ge0$, with the following properties.

(I) $\f(0)=0$, $\f(t)$ is increasing and differentiable as $ t>0$;

(II) $\f(t)\to\infty$ as $t\to\infty$;

(III) $\f'(t),\ t>0$, is monotonic (that is $\f(t)$ is either convex or concave);

(IV) $\f(t)$ satisfies the doubling condition
\begin{equation}\label{f21}
\f(2t)\le2^s\f(t)\ \text{ for all }\ t\ge0\ \text{and for some }s>0\ \text{depending on } \f.
\end{equation}
\end{definition}
Our main result is the following generalization of Theorem~2.7 in \cite{ENV}.

\begin{theorem}\label{th21}
(i) Let $\f\in\Phi$. For any compact set $E\subset\R^d$,

\begin{equation}\label{f22}
\g_{\f,+}(E)\ge c\,\sup\|\mu\|^{3/2}\bigg[\int_{\R^d}W_\f^\mu(x)\,d\mu(x)\bigg]^{-1/2}.
\end{equation}

(ii) Suppose that $\f\in\Phi$, $\f$ is concave, and $s\in(0,1)$, where  $s$ is the exponent in  \eqref{f21}. Then
\begin{equation}\label{f23}
\g_{\f,+}(E)\le C\,\sup\|\mu\|^{3/2}\bigg[\int_{\R^d}W_\f^\mu(x)\,d\mu(x)\bigg]^{-1/2},
\end{equation}
where the supremum is taken over all positive Radon measures supported by $E$, and the constants $c,C$ depend only on $d$, $\f$.
\end{theorem}

\begin{corollary}\label{cor22}
Let a function $\f\in\Phi$ be such that $\varphi(t) =(\log \frac{1}{t})^{-1/2}$, $0 < t\le  e^{-3/2}$,  $\f$ is concave, and $s\in(0,1)$. Then
\begin{equation}\label{f24}
\gamma_{\f,+}(E) \approx C_{\frac{2}{3}d, \frac{3}{2}} (E)
\end{equation}
for every compact set $E$ with diam$(E)\le1$.
\end{corollary}

 Theorem \ref{th21} can be viewed as application of the following results. The next theorem is a generalization of \eqref{E:Wolff_energy}.

\begin{theorem}\label{th23}
(i) Let $\f\in\Phi$, and let $\mu$ be a positive Borel measure (not nesessarily satisfying \eqref{f14}). Then for every measurable set $Q$ in $\R^d$ we have
\begin{equation}\label{f25}
\int_Q |\mathfrak{R}_{\mu,\e}^\f \chi_Q(x)|^2 d\mu(x) \le C\int_Q W_{\f}^{\mu\mid Q}(x) d\mu(x),\quad\e>0,
\end{equation}
where $C$ depends only on $s$.

(ii) Moreover, if $\mu(B(x,t))\to0$ as $t\to0$, $x\in\R^d$, and if  $\f\in\Phi$, $\f$ is concave, and $s\in(0,1)$, then
\begin{equation}\label{f26}
\liminf_{\e\to0}\int_{\R^d} |\mathfrak{R}_{\mu,\e}^\f \mathbf1(x)|^2 d\mu(x) \ge c\int_{\R^d} W_{\f}^{\mu}(x) d\mu(x),\quad c=c(s).
\end{equation}
Both sides of \eqref{f26} might be infinite.
\end{theorem}

We need the notion of Calder\'on-Zygmund (CZ) kernel.
\begin{definition}
Let $\mathcal{X}$ be a metric space. A function $K:\XX\times\XX\to\C$ is said to be a CZ kernel if for some $A>0$ and $\d\in(0,1]$ it satisfies the following two conditions:
\begin{gather}
|K(x,y)|\le\frac{A}{\dist(x,y)^s}\label{f27},\\
|K(x,y)-K(x',y)|,\ |K(y,x)-K(y,x')|\le A\frac{\dist(x,x')^\d}{\dist(x,y)^{s+\d}}\label{f28}
\end{gather}
whenever $x,x',y\in\XX$ and $\dist(x,x')\le\frac12\dist(x,y)$.
\end{definition}
To derive from  Theorem \ref{th23} estimates for norms of the operators $\mathfrak{R}_{\mu,\e}^\f, \ \e>0$, we need a certain version of $T1$ theorem. Our kernel $\frac{y-x}{|y-x|}\,\frac1{\f(|y-x|)}$ is not a CZ kernel in $\R^d$ with the Euclidean distance consistent with the condition $\mu(B(x,r))\le r^s$, and we can not use \cite{NTV98}, \cite{NTV2002}. However one can obtain the desired theorem verifying that arguments in \cite{ENV} and \cite{NTV98} work not only for the Riesz kernel but for our generalized kernel as well. Professor F.~Nazarov suggested another approach. He observed that the set $\R^d$ endowed with the distance $\dist(x,y)=\psi(|x-y|)$ defined below, is a metric space, and our kernel is a CZ kernel in this space (see Lemma \ref{le32}). Here $\psi(r)=\inf\sum_i\f(r_i)^{1/s}$, were the infimum is taken over all finite sequences $\{r_i\}$, $r_i>0$, such that $\sum_i r_i=r$. We are grateful to Fedor Nazarov for this suggestion and for the permission to use it in our paper.

Realizing this idea, we obtain the first part of Theorem \ref{th23} as a particular case of the more general result -- Theorem \ref{th31}. Then we prove the following (weakened) version of $T1$ theorem for CZ operators in metric spaces. As before, we denote by $\Sigma_s$ the class of finite nonnegative Borel measures $\eta$ in a metric space $\XX$ such that $\eta(\BB(x,r))\le r^s$, $x\in\XX,\ r>0$, were $\BB(x,r)=\{y\in\XX:\dist(x,y)<r\}$.

\begin{theorem}\label{th24}
Let $\XX$ be a separable metric space and let $\eta\in\Sigma_s$. Set 
$$
\mathfrak{R}_{\eta,\e}^Kf(x)=\int_{\XX\setminus\BB(x,\e)}K(x,y)f(y)\,d\eta(y),
$$
where $K(x,y)$ is a CZ kernel with the same parameter $s$ in \eqref{f28}. Suppose that
\begin{equation}\label{f29}
\|\mathfrak{R}_{\eta,\e}^K\chi_Q\|_{L^2(\eta\mid Q)}^2\le C\eta(Q),\quad C=C(A,s,\d),
\end{equation}
for every measurable set $Q$. Then the operators $\mathfrak{R}_{\eta,\e}^K$ are uniformly bounded with respect to $\e$, that is
\begin{equation}\label{f210}
\|\mathfrak{R}_{\eta,\e}^K\|_{L^2(\eta)\to L^2(\eta)}\le C',\quad C'=C'(A,s,\d).
\end{equation}
\end{theorem}
In the spaces of homogeneous type, an even the better result is known. In particular, one may assume \eqref{f29} only for cubes or balls. It is the famous $T1$ theorem of David-Journ\'e  (see \cite{DJ} for the Euclidean setting and \cite{C} for homogeneous setting). The nonhomogeneous setting was treated by Nazarov, Treil and Volberg in \cite{NTV97} and \cite{NTV2002}, but only for the Euclidean case. More general kernels in $\R^d$ were considered in \cite{HM10}, and we might use this result to prove Theorem \ref{th21}. But we prefer another approach based on Theorem \ref{th24}. In spite of the references to \cite{ENV} and \cite{NTV98} in our proof of Theorem \ref{th24}, we believe that this proof is still shorter than the proof in \cite{HM10}. Unlike the result in \cite{HM10}, our Theorem \ref{th24} covers far more than the Euclidean case. Note that Theorem \ref{th24} and Theorem 2.1 in \cite{HM10} do not imply each other. 

In the sequel we assume that $K$ antisymmetric, that is $K(x,y)=-K(y,x)$.

Using Theorem \ref{th24}, we obtain a generalization of Theorem 2.6 in \cite{ENV}.

\begin{theorem}\label{th25} (i) Let $\mu$ be a Borel measure in a separable metric space $\XX$, and let $K(x,y)$ be an antisymmetric CZ kernel. Then
\begin{equation}\label{f211}
\pmb|\mathfrak{R}_{\mu}^K\pmb|^2\le C\sup_{x\in\supp\mu}W_s^\mu(x),\quad
W_s^\mu(x)=\int_0^\infty\bigg[\frac{\mu(\BB(x,r))}{r^s}\bigg]^2\,\frac{dr}{r},
\end{equation}
where $C$ depends only on the parameters $A,s,\d$ of a kernel $K(x,y)$.

(ii) On the other hand, if $\mu(B(x,t))\to0$ as $t\to0$, $x\in\R^d$, and if  $\f\in\Phi$, $\f$ is concave, and $s\in(0,1)$, then
\begin{equation}\label{f212}
\pmb|\mathfrak{R}_{\mu}^\f\pmb|^2\ge\frac{c}{\|\mu\|}\int_{\R^d}W_\f^\mu(x)\,d\mu(x),
\quad c=c(s).
\end{equation}
\end{theorem}
We remark that the Wolff potential in general metric spaces has the same form as in the Euclidean space. The preceding results allow us to estimate the so-called operational capacity $\g_{K,op}$ and the capacity $\g_{K,\ast}$ defined by the equalities
\begin{align}
\g_{K,op}(E)&:=\sup\{\|\mu\|:\mu\in\Sigma_s,\ \supp\mu\subset E,\ \pmb|\mathfrak{R}_{\mu}^K\pmb|\le1\},\label{f213}\\
\g_{K,\ast}(E)&:=\sup\{\|\mu\|:\mu\in\Sigma_s,\ \supp\mu\subset E,\ \mathfrak{R}_{\mu,\ast}^K{\bf1}(x)\le1,\ x\in\XX\},\label{f214}
\end{align}
where $\mathfrak{R}_{\mu,\ast}^K{\bf1}(x):=\sup_{\e>0}|\mathfrak{R}_{\mu,\e}^K{\bf1}(x)|$. In the case $\XX=\R^d$ with the distance $\dist(x,y)=\psi(|y-x|)$, and 
$$
K(x,y)=K_\f(x,y)=\frac{y-x}{|y-x|}\,\frac{1}{\f(|y-x|)},
$$
we write $\f$ instead of $K$: $\g_{\f,op}$ and so on. The theorems below establish connections between these capacities.

Following \cite{HM09} we say that a metric space is geometrically doubling if every open ball $\BB(x,r)$ can be covered by at most $N$ balls of radius $r/2$, where $N<\infty$ is independent of $x,r$.

\begin{theorem}\label{th26} Let $\XX$ be a compact Hausdorff geometrically doubling metric space, and let $K$ be an antisymmetric CZ kernel. Then for every bounded closed set $E\subset\XX$,
\begin{equation}\label{f215}
\g_{K,op}(E)\approx\g_{K,\ast}(E),
\end{equation}
where the constants of comparison depend only on the parameters of $K$ and on $N$.
\end{theorem}
The related result in $\R^d$ for $\f(t)=t^s$ was obtained by Volberg \cite[Chapter 5]{V}.

\begin{theorem}\label{th27} Suppose that $\f\in\Phi$, and there is $\L=\L(\f)>0$ for which
\begin{equation}\label{f216}
\int_0^r\frac{t^{d-1}\,dt}{\f(t)}<\L\frac{r^d}{\f(r)},\quad r>0.
\end{equation}
Then
\begin{equation}\label{f217}
\g_{\f,\ast}(E)\le\g_{\f,+}(E)\le C\g_{\f,\ast}(E)
\end{equation}
with $C$ depending only on $s$ and $\L$.
\end{theorem}
For example, if $\f(t)=t^s$, then \eqref{f216} means that $s<d$. A certain relation between $\f$ and $d$ is natural, because in the case $\liminf_{t\to0}t^{-d}\f(t)=0$, the class $\Sigma_\f$ consists of only zero measure.

We prove Theorem \ref{th23} in Section 3, and Theorems \ref{th24}, \ref{th25} in Section 4. Theorems \ref{th26}, \ref{th27} are proved in Section 5. The concluding Section 6 contains proofs of Theorem \ref{th21} and of Corollary \ref{cor22}.

\section{Proof of Theorem \ref{th23} and related results}%Section 3

\begin{theorem}\label{th31} Let $\mu$ be a positive Borel measure in a metric space $\XX$, and let $K(x,y)$ be an antisymmetric CZ kernel. Then for every measurable set $Q$ in $\XX$ we have
\begin{equation}\label{f31}
\int_Q |\mathfrak{R}_{\mu,\e}^K \chi_Q(x)|^2 d\mu(x) \le C\int_Q W_{s}^{\mu\mid Q}(x)\,d\mu(x),\quad\e>0,\end{equation}
where the Wolff potential is defined in \eqref{f211}, and $C$ depends only on the parameters $A,s,\d$ of a kernel $K(x,y)$.
\end{theorem}

\begin{proof}
Our arguments are similar to those in the proof of Theorem 2.6 in \cite{ENV}, but there are essential differences as well. For estimation of the right hand side of \eqref{f31} the measure on $\XX\setminus Q$ is unessential. Thus, we may assume that $\mu$ is concentrated on $Q$, and write $\mu$ instead of $\mu|Q$. Also without loss of generality we may assume that
\begin{equation}\label{f32}
\int_0^\infty\bigg[\frac{\mu(\BB(x,r))}{r^s}\bigg]^2\,\frac{dr}{r}<\infty\quad\mu\text{-a.~e.}
\end{equation}
Otherwise (\ref{f31}) becomes trivial.

Let $\e>0$ and a measurable set $Q$ be given. We set
$$\aligned
\UU&=\{(x,y,z)\in Q^3:\dist(x,y)>\e,\ \dist(x,z)>\e\},\\
\UU_1&=\{(x,y,z)\in Q^3:\dist(x,y)\ge\dist(x,z)>\e\},\\
\UU_2&=\{(x,y,z)\in Q^3:\e<\dist(x,y)<\dist(x,z)\},\\
\UU_{1,1}&=\{(x,y,z)\in Q^3:\dist(x,y)\ge\dist(x,z)>\e,\ \dist(y,z)\ge\dist(x,z)\},\\
\UU_{1,2}&=\{(x,y,z)\in Q^3:\dist(x,y)\ge\dist(x,z)>\e,\ \dist(y,z)<\dist(x,z)\}.
\endaligned$$
Then
$$\aligned
\int_{Q}|\mathfrak{R}_{\mu,\e}^K\chi_Q(x)|^2\,d\mu(x)&=\iiint_{\UU}
K(x,y)K(x,z)\,d\mu(z)\,d\mu(y)\,d\mu(x)\\
&=\iiint_{\UU_1}+\iiint_{\UU_2}=:I_1+I_2.
\endaligned$$
Estimates for $I_1,I_2$ are analogous. It is enough to estimate $I_1$. We have
$$\aligned
|I_1|\le&\biggl|\iiint_{\UU_{1,1}}K(x,y)K(x,z)\,d\mu(z)\,d\mu(y)\,d\mu(x)\biggr|\\
+&\biggl|\iiint_{\UU_{1,2}}K(x,y)K(x,z)\,d\mu(z)\,d\mu(y)\,d\mu(x)\biggr|=:I_{1,1}+I_{1,2}.
\endaligned$$
We put the absolute value in $I_{1,2}$ inside the integral. Since $\dist(x,z)>\frac12\dist(x,y)$ in $I_{1,2}$, \eqref{f27} yields the estimate
\begin{equation}\label{f33}\aligned
I_{1,2}&\le\int_Q\int_{\dist(x,y)>0}\frac{2^sA^2}{\dist(x,y)^{2s}}\,\mu(\BB(x,\dist(x,y)))\,d\mu(y)\,d\mu(x)\\
&=2^sA^2\int_Q\int_0^\infty\frac1{r^{2s}}\,\mu(\BB(x,r))\,d\mu(\BB(x,r))\,d\mu(x)\\
&=2^sA^2\int_Q\int_0^\infty\frac1{r^{2s}}\,d\biggl[\frac{\mu(\BB(x,r))^2}2\biggr]\,d\mu(x).
\endaligned\end{equation}
From (\ref{f32}) one can easily deduce that
\begin{equation}\label{f34}
\lim_{r\to0}\frac{\mu(\BB(x,r))}{r^{s}}=0,\quad
\lim_{r\to\infty}\frac{\mu(\BB(x,r))}{r^{s}}=0,\quad\mu\text{-a.~e.}
\end{equation}
Integrating by parts in the last integral of (\ref{f33}) we get
$$
I_{1,2}\le s2^sA^2\int_Q\int_0^\infty\bigg[\frac{\mu(\BB(x,r))}{r^s}\bigg]^2\,\frac{dr}{r}\,d\mu(x)=s2^sA^2\int_QW_{s}^{\mu\mid Q}(x)\,d\mu(x).
$$

Let us estimate $I_{1,1}$. By the symmetry of $\UU_{1,1}$ with respect to $z,x$ we have
$$\aligned
I_{1,1}&=\frac12\biggl|\iiint_{\UU_{1,1}}(K(x,y)K(x,z)+K(z,y)K(z,x))\,d\mu(z)\,d\mu(y)\,d\mu(x)\biggr|\\
&\le\frac12\iiint_{\UU_{1,1}}|K(x,z)|\cdot|K(x,y)-K(z,y)|\,d\mu(z)\,d\mu(y)\,d\mu(x)
\endaligned$$
(recall that $K(x,y)$ is antisymmetric). If $\dist(x,z)\le\frac12\dist(x,y)$, then from \eqref{f27}, \eqref{f28} we deduce
$$
|K(x,z)|\cdot|K(x,y)-K(z,y)|
\le\frac{A}{\dist(x,z)^s}\cdot \frac{A\dist(x,z)^\d}{\dist(x,y)^{s+\d}}
=\frac{A^2}{\dist(x,y)^{s+\d}\dist(x,z)^{s-\d}}.
$$
 If $\dist(x,z)>\frac12\dist(x,y)$, then we derive the analogous estimate directly from \eqref{f27} with another constant $C=C(s)$ instead of $A^2$. Hence,
$$\aligned
I_{1,1}&\le C\iiint_{\UU_{1,1}}\frac{1}{\dist(x,y)^{s+\d}\dist(x,z)^{s-\d}}
\,d\mu(z)\,d\mu(y)\,d\mu(x)\\
&\le C\int_Q\int_{\dist(x,y)\ge\e}\frac1{\dist(x,y)^{s+\d}}
\biggl[\int_\e^{\dist(x,y)}\frac{d\mu(\BB(x,t))}{t^{s-\d}}\biggr]\,d\mu(y)\,d\mu(x)\\
&\le C\int_Q\int_0^\infty\frac1{r^{s+\d}}\biggl[\int_0^{r}\frac{d\mu(\BB(x,t))}{t^{s-\d}}\biggr]\,d\mu(\BB(x,r))\,d\mu(x),\quad C=C(A,s).
\endaligned$$
Set
$$
H_x(r):=\int_0^{r}\frac{d\mu(\BB(x,t))}{t^{s-\d}}.
$$
Then the last expression can be written in the form
\begin{equation}\label{f35}
C\int_Q\int_0^\infty\frac1{r^{2\d}}H_x(r)\,dH_x(r)\,d\mu(x)
=\frac{C}2\int_Q\int_0^\infty\frac{dH_x(r)^2}{r^{2\d}}\,d\mu(x).
\end{equation}
Obviously,
\begin{equation}\label{f36}
H_x(r)=\frac{\mu(\BB(x,r))}{r^{s-\d}}+(s-\d)\int_0^{r}\frac{\mu(\BB(x,t))}{t^{s-\d+1}}\,dt,
\end{equation}
and by (\ref{f34}) we have
$$
\lim_{r\to\infty}\frac{H_x(r)}{r^\d}=0,\quad \lim_{r\to0}\frac{H_x(r)}{r^\d}=0\quad \mu\text{-a.~e.}
$$
Thus,
\begin{equation}\label{f37}\aligned
\int_0^\infty\frac{dH_x(r)^2}{r^{2\d}}&=2\d\int_0^\infty\frac{H_x(r)^2}{r^{2\d+1}}\,dr
\overset{(\ref{f36})}\le4\d\int_0^\infty\bigg[\frac{\mu(\BB(x,r))}{r^s}\bigg]^2\,\frac{dr}{r}\\
&+4\d(s-\d)^2\int_0^\infty\frac1{r^{2\d+1}}\biggl[\int_0^{r}\frac{\mu(\BB(x,t))}{t^{s-\d+1}}\,dt\biggr]^2dr.
\endaligned\end{equation}
The first term in the right hand side of (\ref{f37}) is what we need. Let us estimate the second term.
By the Cauchy--Bunyakovskii--Schwarz inequality,
$$
\biggl[\int_0^{r}\frac{\mu(\BB(x,t))}{t^{s-\d+1}}\,dt\biggr]^2\le
\int_0^{r}\bigg[\frac{\mu(\BB(x,t))}{t^{s+\frac{1-\d}2}}\bigg]^2\,dt\cdot\int_0^r\frac{dt}{t^{1-\d}}=\frac{r^\d}{\d}\int_0^{r}\bigg[\frac{\mu(\BB(x,t))}{t^{s}}\bigg]^2\frac{dt}{t^{1-\d}}.
$$
Applying integration by parts, we obtain the estimate
\begin{multline*}
\int_0^\infty\frac1{r^{2\d+1}}\biggl[\int_0^{r}\frac{\mu(\BB(x,t))}{t^{s-\d+1}}\,dt\biggr]^2dr
\le\frac1{\d}\int_0^\infty\biggl\{\int_0^{r}\bigg[\frac{\mu(\BB(x,t))}{t^s}\bigg]^2\frac{dt}{t^{1-\d}}\biggr\}\frac{dr}{r^{1+\d}}\\
=-\frac1{\d^2}\biggl(\frac1{r^\d}\int_0^{r}\bigg[\frac{\mu(\BB(x,t))}{t^s}\bigg]^2\frac{dt}{t^{1-\d}}\biggr)\bigg|_0^\infty
+\frac1{\d^2}\int_0^{\infty}\bigg[\frac{\mu(\BB(x,r))}{r^s}\bigg]^2\,\frac{dr}r.
\end{multline*}
According to (\ref{f34}), the substitution of limits gives zero. Thus, (see (\ref{f37}))
$$
\int_0^\infty\frac{dH_x(r)^2}{r^{2\d}}<C(s,\d)\int_0^{\infty}\bigg[\frac{\mu(\BB(x,r))}{r^s}\bigg]^2\,\frac{dr}r.
$$
Now (\ref{f35}) yields (\ref{f31}), and Theorem \ref{th31} is proved.
\end{proof}

\begin{lemma}\label{le32}
Suppose that $\f\in\Phi$. Let
 $$
\psi(r)=\inf\sum_i\f(r_i)^{1/s},
$$
were $s$ is the exponent in \eqref{f21}, and  the infimum is taken over all finite sequences $\{r_i\}$, $r_i>0$, such that $\sum_i r_i=r$. The following statements hold.

(i) The set of points $x\in\R^d,\ d\ge1$, with the distance $\dist(x,y)=\psi(|x-y|)$, is a metric space.

(ii) The kernel
$$
K_\f(x,y)=\frac{y-x}{|y-x|}\,\frac1{\f(|y-x|)}
$$
is a CZ kernel in the metric space $\XX$ defined above in (i) with the parameters $s$ from \eqref{f21}, $\d=1$, and $A=A(s)$.

(iii) The condition $\mu(B(x,r))\le\f(r)$ implies that $\mu(\BB(x,r))\le C(s)r^s,\ r>0$, where $B(x,r)$ is a Euclidean ball, and $\BB(x,r)$ is a ball in $\XX$. Conversely, if $\mu(\BB(x,r))\le r^s,\ r>0$, then $\mu(B(x,r))\le\f(r)$.
\end{lemma}

\begin{proof}
(i) We prove that
\begin{equation}\label{f38}
\tfrac12\f(r)^{1/s}\le\psi(r)\le\f(r)^{1/s},\quad r\ge0.
\end{equation}
Given $r>0$, find $\overline{r}\in(0,r]$ for which
$$
\min_{0<t\le r}\frac{\f(t)^{1/s}}{t}=\frac{\f(\overline{r})^{1/s}}{\overline{r}}.
$$
The minimum is attained for some $\overline{r}\in(r/2,r]$, because by \eqref{f21} we have
\begin{equation}\label{f39}
\frac{\f(2t)^{1/s}}{2t}\le\frac{\f(t)^{1/s}}{t},\quad t>0.
\end{equation}
For every sequence $\{r_i\}$ with $\sum_i r_i=r$ we get
$$
\sum_i\f(r_i)^{1/s}=\sum_i\frac{\f(r_i)^{1/s}}{r_i}\,r_i
\ge\frac{\f(\overline{r})^{1/s}}{\overline{r}}\,r\overset{(\ref{f39})}\ge
\frac{\f(2\overline{r})^{1/s}}{2\overline{r}}\,r>\frac{\f(r)^{1/s}}{2r}\,r,
$$
and the first inequality in \eqref{f38} is proved. The second one is trivial.

Relations \eqref{f38} imply that $\dist(x,y)=0\iff x=y$. The definition of $\psi$ yields the property $\psi(a+b)\le\psi(a)+\psi(b)$, $a,b>0$. Hence,
$$
\dist(x,z)=\psi(|x-z|)\le\psi(|x-y|+|y-z|)\le\dist(x,y)+\dist(y,z),
$$
and the part (i) is proved.

(ii) The property \eqref{f27} with $A=1$ easily follows from \eqref{f38}:
$$
|K_\f(x,y)|=\frac1{\f(|y-x|)}\le\frac1{\psi(|y-x|)^s}=\frac1{\dist(x,y)^s}.
$$
To establish \eqref{f28} we need the following property of $\f$:
\begin{equation}\label{f310}
\frac{t_1}{t_2}<\frac{2\f(t_1)^{1/s}}{\f(t_2)^{1/s}},\quad 0<t_1\le t_2.
\end{equation}
Indeed, take the integer $k\ge0$ for which $2^kt_1\le t_2<2^{k+1}t_1$. Then
$$
\frac{\f(t_2)^{1/s}}{t_2}<2\,\frac{\f(2^{k+1}t_1)^{1/s}}{2^{k+1}t_1}\overset{(\ref{f39})}\le
\frac{2\f(t_1)^{1/s}}{t_1}.
$$
Let $x,x',y$ be such that $\dist(x,x')\le\frac12\dist(x,y)$. Set $a=x-y,\ b=x'-y$. We have
$$\aligned
|K_\f(x,y)-K_\f(x',y)|&=\biggl|\frac{a}{|a|}\,\frac1{\f(|a|)}-\frac{b}{|b|}\,\frac1{\f(|b|)}\biggr|\\
&\le\biggl|\frac{a}{|a|}-\frac{b}{|b|}\biggr|\,\frac1{\f(|a|)}+\biggr|\frac1{\f(|a|)}-\frac1{\f(|b|)}\biggr|\\
&\le\biggr[\frac{|a-b|}{|a|}+|b|\biggl|\frac1{|a|}-\frac1{|b|}\biggr|\biggr]\frac1{\f(|a|)}+
\frac{|\f(|b|)-\f(|a|)|}{\f(|a|)\f(|b|)}\\
&\le\frac{2|a-b|}{|a|}\frac1{\f(|a|)}+\frac{\f'(\xi)|a-b|}{\f(|a|)\f(|b|)},
\endaligned$$
where $\xi$ is a number between $|a|$ and $|b|$. Suppose that $|a|\le|b|$. If $\f'(t)$ is nonincreasing then $\f'(\xi)\le\f'(|a|)\le\f(|a|)/|a|$. Hence,
$$
|K_\f(x,y)-K_\f(x',y)|\le\frac{|a-b|}{|a|}\biggl[\frac2{\f(|a|)}+\frac1{\f(|b|)}\biggr]
\le\frac{3|a-b|}{|a|\f(|a|)}.
$$
 If $\f'(t)$ is nondecreasing, we have
$$
\f'(\xi)\le\f'(|b|)\le\frac{\f(2|b|)-\f(|b|)}{|b|}\overset{(\ref{f21})}<\frac{2^s\f(|b|)}{|a|}.
$$
In this case
$$
|K_\f(x,y)-K_\f(x',y)|\le(2+2^s)\frac{|a-b|}{|a|\f(|a|)}.
$$
From (\ref{f310}) and (\ref{f38}) we get
$$
\frac{|a-b|}{|a|\f(|a|)}<\frac{2\f(|a-b|)^{1/s}}{\f(|a|)^{1/s}}\,\frac1{\f(|a|)}
<\frac{4\psi(|a-b|)}{\psi(|a|)^{s+1}}=\frac{4\dist(x,x')}{\dist(x,y)^{s+1}}.
$$
We consider the case $|b|\le|a|$ in the same way, taking into account that
$\dist(x',y)\le\dist(x,y)\le2\dist(x',y)$.

(iii) The last statements follow from the obvious relation $B(x,r)=\BB(x,\psi(r))$. Let $t=\psi(r)$. Then
$\mu(\BB(x,t))=\mu(B(x,r))\le\f(r)\overset{(\ref{f38})}\le2^s\psi(r)^s=2^st^s$. Conversely,
$\mu(B(x,r))=\mu(\BB(x,t))\le t^s=\psi(r)^s\overset{(\ref{f38})}\le\f(r)$.
Lemma \ref{le32} is proved.
\end{proof}

We consider the quantity
\begin{multline}\label{f311}
2p_\f(x_1,x_2,x_3):=\\
\sum_\sigma\frac{x_{\sigma(2)}-x_{\sigma(1)}}{|x_{\sigma(2)}-x_{\sigma(1)}|}
\cdot\frac1{\f(|x_{\sigma(2)}-x_{\sigma(1)}|)}\cdot
\frac{x_{\sigma(3)}-x_{\sigma(1)}}{|x_{\sigma(3)}-x_{\sigma(1)}|}
\cdot\frac1{\f(|x_{\sigma(3)}-x_{\sigma(1)}|)}\,,
\end{multline}
where $x_1,x_2,x_3$ are given three distinct points in $\R^d$ and the sum is taken over the six permutations of the set $\{1,2,3\}$. This quantity is an analog of Menger curvature \cite{M}. It was observed in \cite{P}, p.~952, that one can define $p_\f(x_1,x_2,x_3)$ as the sum in (\ref{f311}) taken over only the three permutations (1,2,3), (2,3,1) and (3,1,2), since the other three permutations give the same terms in (\ref{f31}). Later on we also will write $x,y,z$ instead of $x_1,x_2,x_3$.

\begin{lemma}\label{le33}
Let $x,y,z$ be three distinct points in $\R^d$, $d\ge 1$, and let $\f(t)$, $t\ge0$, be an increasing function with $\f(0)=0$. Set $a=|y-x|$, $b=|z-y|$ and $c=|z-x|$. If $a\ge b\ge c$, then
\begin{equation}\label{f312}
p_\f(x,y,z)\le\frac{c}{b\f(b)\f(c)}+\frac{2}{\f(a)\f(b)}\,.
\end{equation}

Moreover, if $\f(t)\in\Phi$, $\f$ is concave, and $s\in(0,1)$, then
\begin{equation}\label{f313}
p_\f(x,y,z)>\frac{1-2^{s-1}}{4\f(a)\f(b)}\,.
\end{equation}
\end{lemma}
One can derive \eqref{f25} directly from \eqref{f312}. But we have Theorem \ref{th31} and Lemma \ref{le32}. Thus, we will not use \eqref{f312} in the sequel, and give the short proof of this inequality for completeness.

\begin{proof}
Let $\a,\beta,\g$ be the angles opposite to sides $a,b,c$ respectively. Since
$$
\cos\alpha=\frac{b^2+c^2-a^2}{2bc}\,,\quad\cos\beta=\frac{a^2+c^2-b^2}{2ac}\,,\quad\cos\g=\frac{a^2+b^2-c^2}{2ab}\,,
$$
we have
\begin{equation}\label{f314}
\begin{split}
p_\f(x,y,z)&=[\f(a)\f(b)\f(c)]^{-1}(\f(a)\cos\a+\f(b)\cos\beta+\f(c)\cos\g)\\
&=\frac1{2\f(b)\f(c)}\bigg[\frac{b^2+c^2-a^2}{bc}+\frac{\f(b)}{\f(a)}\,\frac{a^2+c^2-b^2}{ac}
+\frac{\f(c)}{\f(a)}\,\frac{a^2+b^2-c^2}{ab}\bigg].
\end{split}
\end{equation}
Let us prove \eqref{f312}.
$$\aligned
p_\f(x,y,z)&=\frac1{2\f(b)\f(c)}\bigg[c\bigg(\frac1b+\frac{\f(b)}{\f(a)}\,\frac1a\bigg)
+\frac{a^2-b^2}{c}\,\bigg(\frac{\f(b)}{\f(a)}\,\frac1a-\frac1b\bigg)
+\frac{\f(c)}{\f(a)}\,\frac{a^2+b^2-c^2}{ab}\bigg]\\
&\le\frac1{2\f(b)\f(c)}\bigg[c\,\frac2b+\frac{\f(c)}{\f(a)}\,\frac{2a^2}{ab}\bigg]
\le\frac{c}{b\f(b)\f(c)}+\frac2{\f(a)\f(b)}\,,
\endaligned$$
since $a/2\le b\le a$ and
$$
\frac{\f(b)}{\f(a)}\,\frac1a-\frac1b=\frac{b\f(b)-a\f(a)}{ab\f(a)}\le0.
$$

To get the lower bound \eqref{f313}, we set $u=b/a,\ v=c/a$ and write \eqref{f314} as
\begin{equation}\label{f315}
\begin{split}
&p_\f(x,y,z)\\
&=\frac1{2uv\f(b)\f(c)}\bigg[u^2+v^2-1+\frac{\f(b)}{u\f(a)}\,u^2(1+v^2-u^2)
+\frac{\f(c)}{v\f(a)}\,v^2(1+u^2-v^2)\bigg]\\
&=\frac1{2uv\f(b)\f(c)}\bigg[\bigg(\frac{\f(b)}{u\f(a)}-1\bigg)\,u^2(1+v^2-u^2)
+\bigg(\frac{\f(c)}{v\f(a)}-1\bigg)\,v^2(1+u^2-v^2)\\
&+\{u^2+v^2-1+u^2(1+v^2-u^2)+v^2(1+u^2-v^2)\}\bigg].
\end{split}
\end{equation}
The expression in braces is equal to
$$\aligned
D&:=2u^2+2v^2+2u^2v^2-u^4-v^4-1\\
&=4v^2-[(u^2-v^2)^2-2(u^2-v^2)+1]=4v^2-(u^2-v^2-1)^2\\
&=(1+v-u)(1+v+u)(1+u-v)(u+v-1)\ge0.
\endaligned$$
Set $\l=2^{s-1}$ and consider two cases.

{\bf Case 1.} $v\ge(3-\l)/4$. Then
$$
D\ge v\cdot2\cdot1\cdot\frac{1-\l}2=v(1-\l).
$$
Since $\f(t)/t$ is nonincreasing, we have
$$
\frac{\f(b)}{u\f(a)}-1=\bigg(\frac{\f(b)}{b}-\frac{\f(a)}{a}\bigg)\cdot\frac{a}{\f(a)}\ge0,\quad
\frac{\f(c)}{v\f(a)}-1\ge0.
$$
From \eqref{f315} we get
$$
p_\f(x,y,z)\ge\frac{1-\l}{2\f(b)\f(c)}\ge\frac{1-\l}{2\f(a)\f(b)}\,.
$$

{\bf Case 2.} $v<(3-\l)/4$. Then $a>4c/(3-\l)$. Hence,
$$\aligned
&\frac{\f(c)}{v\f(a)}-1=\bigg(\frac{\f(c)}{c}-\frac{\f(a)}{a}\bigg)\,\frac{a}{\f(a)}
\ge\bigg(\frac{\f(c)}{c}-\frac{\f(4c/(3-\l))}{4c/(3-\l)}\bigg)\,\frac{a}{\f(a)}\\
&>\bigg(\frac{\f(c)}{c}-\frac{\f(2c)(3-\l)}{4c}\bigg)\,\frac{a}{\f(a)}
\overset{\eqref{f21}}\ge\bigg(\frac{\f(c)}{c}-\frac{2\l\f(c)(3-\l)}{4c}\bigg)\,\frac{a}{\f(a)}\\
&=\frac{a}{\f(a)}\,\frac{\f(c)}{c}\,\bigg(1-\frac{\l(3-\l)}{2}\bigg)
=\frac{\f(c)}{2v\f(a)}(1-\l)(2-\l)>\frac{\f(c)}{2v\f(a)}(1-\l).
\endaligned$$
Now \eqref{f315} implies the estimate
$$
p_\f(x,y,z)>\frac1{2uv\f(b)\f(c)}\,\frac{\f(c)}{2v\f(a)}(1-\l)v^2\ge\frac{1-\l}{4\f(a)\f(b)}\,.
$$
Lemma \ref{le33} is proved.
\end{proof}

\begin{proof}[Proof of Theorem \ref{th23}]
(i) Obviously, $\mathfrak{R}_{\mu,\e}^\f f(x)=\mathfrak{R}_{\mu,\psi(\e)}^Kf(x)$ with
$K(x,y)=\frac{y-x}{|y-x|}\,\frac{1}{\f(|y-x|)}$. Using the substitution $r=\psi(t)$ and applying \eqref{f31}, we get
\begin{multline*}
\int_Q |\mathfrak{R}_{\mu,\e}^\f \chi_Q(x)|^2 d\mu(x)\le C\int_Q \int_0^\infty\biggl[
\frac{(\mu\!\mid\! Q)(\BB(x,r))}{r^s}\biggr]^2\frac{dr}r\,d\mu(x)\\
=C\int_Q \int_0^\infty\biggl[
\frac{(\mu\!\mid\! Q)(B(x,t))}{\psi(t)^s}\biggr]^2\frac{d\psi(t)}{\psi(t)}\,d\mu(x)
\le2^{2s} C\int_QW_\f^{\mu\mid Q}(x)\,d\mu(x),\quad C=C(s);
\end{multline*}
in the last inequality we use \eqref{f38} and integration by parts.

(ii) We use the standard symmetrization arguments (see \cite{MV}). Fix $\e>0$ and set
$$\aligned
\UU&=\{(x,y,z)\in (\R^d)^3:|y-x|>\e,\ |z-x|>\e\},\\
\UU_1&=\{(x,y,z)\in (\R^d)^3:|y-x|>\e,\ |z-x|>\e,\ |z-y|\le\e\},\\
\O&=\UU\setminus\UU_1=\{(x,y,z)\in (\R^d)^3:|y-x|>\e,\ |z-x|>\e,\ |z-y|>\e\}.
\endaligned$$
Then
\begin{align}
\int_{\R^d}|\mathfrak{R}_{\mu,\e}^\f\mathbf1(x)|^2\,d\mu(x)&=\iiint_{\UU}
\frac{y-x}{|y-x|}\,\frac{1}{\f(|y-x|)}\,\frac{z-x}{|z-x|}\,\frac{1}{\f(|z-x|)}\,d\mu(z)\,d\mu(y)\,d\mu(x)\notag\\
&=\iiint_{\O}+\iiint_{\UU_1}=:I_1+I_2.\label{f316}
\end{align}
The set $\O$ is symmetric with respect to permutations of $x,y,z$. Hence,
\begin{equation}\label{f317}
I_1=\frac13\iiint_{\O}p_\f(x,y,z)\,d\mu(z)\,d\mu(y)\,d\mu(x).
\end{equation}

For $(x,y,z)\in\UU_1$, the angle between the vectors $y-x$ and $z-x$ is acute. Hence, the triple integral over $\UU_1$ in \eqref{f316} is positive. The triple integral in \eqref{f317} is greater than or equal to
$$
\iiint_{\O_2}p_\f(x,y,z)\,d\mu(z)\,d\mu(y)\,d\mu(x),
$$
where
$$
\O_2=\{(x,y,z)\in(\R^d)^3:|y-x|>\e,\ \e<|z-x|<|y-x|,\ |z-y|>\e\}.
$$
For triples $(x,y,z)\in\O_2$, the largest side length of the triangle $x,y,z$ does not exceed $2|y-x|$. According to \eqref{f313}, the last triple integral over $\O_2$ is greater than
\begin{multline}\notag
\frac{1-2^{s-1}}4\int_{\R^d}\int_{|y-x|>\e}\int_{z:(x,y,z)\in\O_2}
\frac{1}{\f(2|y-x|)\f(|y-x|)}\,d\mu(z)\,d\mu(y)\,d\mu(x)\\
\overset{\eqref{f21}}\ge\frac{1-2^{s-1}}{4\cdot2^s}\int_{\R^d}\int_{|y-x|>\e}
\frac{\mu(B(x,|y-x|)\setminus\overline{B}(x,\e)\setminus\overline{B}(y,\e))}{\f^2(|y-x|)}\,d\mu(y)\,d\mu(x).
\end{multline}
Set
\begin{equation}\notag
\xi(\e,x,y)=\left\{\begin{array}{ll}
\dfrac{\mu(B(x,|y-x|)\setminus\overline{B}(x,\e)\setminus\overline{B}(y,\e))}{\f^2(|y-x|)},&|y-x|>\e,\\
0,&|y-x|\le\e.
\end{array}\right.
\end{equation}
Clearly,
$$
\xi(\e,x,y)\nearrow\frac{\mu(B(x,|y-x|))}{\f^2(|y-x|)}=:\xi(0,x,y)\text{ as }\e\to0,\ \ |y-x|>0.
$$
 Hence,
$$
\zeta(\e,x):=\int_{|y-x|>\e}\xi(\e,x,y)\,d\mu(y)=\int_{\R^d}\xi(\e,x,y)\,d\mu(y)\nearrow \int_{\R^d}\xi(0,x,y)\,d\mu(y)\ \text{ as }\ \e\to0.
$$
Note that
$$
\int_{\R^d}\xi(0,x,y)\,d\mu(y)=\int_0^\infty\frac{\mu(B(x,t))\,d\mu(B(x,t))}{\f^2(t)}
=W_\f^\mu(x),\ \ x\in\R^d.
$$
The first equality is obvious. If $\mu(B(x,t))/\f(t)\not\to0$ as $t\to0$ or/and $t\to\infty$, then both parts of the last equality are infinite. If $\mu(B(x,t))/\f(t)\to0$ as $t\to0$ and $t\to\infty$, then we obtain this equality integrating by parts. We conclude that
$$
\int_{\R^d}\zeta(\e,x)\,d\mu(x)\nearrow \int_{\R^d}W_\f^\mu(x)\,d\mu(x),
$$
and the proof of Theorem \ref{th23} is completed.
\end{proof}

\section{Proof of Theorems \ref{th24} and \ref{th25}}%Section 4

To avoid the revision of the proof of the nonhomogeneous $T1$ and $Tb$ theorems given by Nazarov, Treil and Volberg in \cite{NTV97}, \cite{NTV2002}, as well as their generalization \cite{HM10}, we will follow the arguments from \cite{ENV}, namely the second approach to Theorem 2.6 in \cite{ENV}. But we will prove a weaker assertion than $T1$ theorem: we assume that \eqref{f29} holds for all measurable sets $Q$, not only for cubes.

\begin{proof}[Proof of Theorem \ref{th24}] The first step is the passage from the trancated operators $\mathfrak{R}_{\eta,\e}^K$ (which are not operators with CZ kernels) to similar operators, but with CZ kernels. Let $\phi(t),\ t\ge0$, be a $C^\infty$ function such that $\phi(t)=0$ as $0\le t\le1$, $\phi(t)=1$ as $t\ge2$, and $0\le\phi'(t)\le2,\ t>0$. Let $\phi_\e(t):=\phi(\frac{t}{\e})$. We prove that if $K(x,y)$ is a CZ kernel with constants $A,s,\d$, then
$$
K_\e(x,y):=\phi_\e(\dist(x,y))K(x,y)
$$
is a CZ kernel with constants $9A,s,\d$. Indeed, the validity of \eqref{f27} for $K_\e$ with the same constant $A$ is obvious. To prove that
\begin{equation}\label{f41}
|K_\e(x,y)-K_\e(x',y)|\le9A\frac{\dist(x,x')^\d}{\dist(x,y)^{s+\d}}
\end{equation}
whenever $\dist(x,x')\le\frac12\dist(x,y)$, we may assume that for at least one point $x$ or $x'$ (say, for $x'$), $\dist(x',y)<2\e$ (otherwise \eqref{f41} follows from \eqref{f28}). Then
$$
\dist(x,y)\le\dist(x',y)+\dist(x',x)\le2\e+\tfrac12\dist(x,y).
$$
Hence, $\dist(x,y)\le4\e$. We have
$$\aligned
&|K_\e(x,y)-K_\e(x',y)|\\
&\le|\phi_\e(\dist(x,y))-\phi_\e(\dist(x',y))|\cdot|K(x,y)|+|\phi_\e(\dist(x',y))|\cdot|K(x,y)-K(x',y)|\\
&\le\frac2{\e} |\dist(x,y)-\dist(x',y)|\cdot|K(x,y)|+A\frac{\dist(x,x')^\d}{\dist(x,y)^{s+\d}}\\
&\le\frac2{\e}\dist(x,x')|K(x,y)|+ A\frac{\dist(x,x')^\d}{\dist(x,y)^{s+\d}}
\le\frac{2A\dist(x,x')}{\e\dist(x,y)^{s}}+A\frac{\dist(x,x')^\d}{\dist(x,y)^{s+\d}}\\
&\le\frac{A}{\dist(x,y)^{s}}\,\bigg[\frac{8\dist(x,x')}{\dist(x,y)}+
\bigg(\frac{\dist(x,x')}{\dist(x,y)}\bigg)^\d\bigg]\le9A\frac{\dist(x,x')^\d}{\dist(x,y)^{s+\d}},
\endaligned$$
and we get \eqref{f41}. In the case $\dist(x,y)\le2\e$ we have $5A$ instead of $9A$. The proof of the analogous estimate for $|K_\e(y,x)-K_\e(y,x')|$ is essentially the same. Set
$$\aligned
\widetilde{\mathfrak{R}}_{\mu,\e}^K f(x)&=\int K_{\e}(x,y)f(y)\,d\mu(y),\quad f\in L^2(\mu),\quad \e>0,\\
\widetilde{\mathfrak{R}}_{\mu,\ast}^K f(x)&=\sup_{\e>0}|\widetilde{\mathfrak{R}}_{\mu,\e}^K f(x)|.
\endaligned$$
We denote by $\widetilde{R}_{\nu,\e}^K$ the corresponding modified $s$-Riesz transform of a finite Borel measure $\nu$:
$$
\widetilde{R}_{\nu,\e}^K(x)=\int K_{\e}(x,y)\,d\nu(y).
$$

The rest of the proof is the same as in \cite{ENV}, starting from inequality (3.11) in \cite{ENV} until the end of proof of Lemma 3.4 with the following minor corrections. All constants $C$ and $C_0$ now depend on CZ constants of the kernel (instead of $d,s$ in \cite{ENV}). The reference [21, Lemma 2.1] after equality (3.12) in \cite{ENV} should be replaced with [21, Lemma 3.1] (in fact, this is a misprint). Finally, the only place where the specific character of the Euclidean metric and of the Riesz kernel is used, is the following simple statement in the proof of Lemma 3.4 in \cite{ENV}. For given $f\in L^1(\eta)$, $\e>0$, $t>0$, one can approximate $f\,d\eta$ by a measure $\nu$ of the form $\nu=\sum_{j=1}^M \a_j\d_{y_j}$, $M\in\mathbb N_+$, $\a_j\in\R$, in such a way that $\|\nu\|\le\|f\|_{L^1(\eta)}$, and
\begin{equation}\label{f42}
\{|\widetilde{\mathfrak{R}}_{\eta,\e}^K f(x)|>t\}\subset\{|\widetilde{R}_{\nu,\e}^K(x)|>\tfrac12 t\}.
\end{equation}
One can easily prove this assertion without the notion of cubes (we do not have it in metric spaces in general), and without equicontinuity. With this end in view we choose $\e'\in(0,\e/4)$ in the following way:
\begin{equation}\label{f43}
| K_{\e}(x,y)- K_{\e}(x,y')|<\frac{t}{4\|f\|_{L^1(\eta)}}
\end{equation}
whenever $\dist(y,y')<\e'$, $x\in\XX$. It is possible because
$$
| K_{\e}(x,y)- K_{\e}(x,y')|<\frac{9A\dist(y,y')^\d}{\dist(x,y)^{s+\d}}
<\frac{A'\dist(y,y')^\d}{\e^{s+\d}},\quad\dist(y,y')<\frac{\e}4
$$
(we recall that $ K_{\e}(x,y)=0$ as $\dist(x,y)<\e$). Let $\{y_i\}$ be a countable everywhere dense subset of $\XX$. Obviously, $\bigcup_i\BB(y_i,\e')=\XX$. Set $Q_1=\BB(y_1,\e')$,
$Q_k=\BB(y_k,\e')\setminus\bigcup_{i=1}^{k-1}\BB(y_i,\e'),\ k=2,\dots$ (possibly, $Q_k=\varnothing$), and take $\a_k=\int_{Q_k}f\,d\eta$. Then $Q_i\cap Q_j\ne\varnothing,\ i\ne j$, and $\bigcup_iQ_i=\XX$. Using \eqref{f43} we have
$$
\sum_{i=1}^\infty\int_{Q_i}| K_{\e}(x,y)- K_{\e}(x,y_i)|\,|f(y)|\,d\eta(y)
<\frac{t}{4\|f\|_{L^1(\eta)}}\,\|f\|_{L^1(\eta)}=\frac t4.
$$
There is $M\in\mathbb N_+$ such that
$$
\sum_{i=M+1}^\infty\int_{Q_i}|K_{\e}(x,y_i)|\,|f(y)|\,d\eta(y)
<\frac{A}{\e^s}\sum_{i=M+1}^\infty\int_{Q_i}|f(y)|\,d\eta(y)<\frac t4,\quad x\in\XX.
$$
Thus,
$$
|\widetilde{\mathfrak{R}}_{\eta,\e}^K f(x)-\widetilde{R}_{\nu,\e}^K(x)|<\frac t4+\frac t4=\frac t2,
$$
and we obtain \eqref{f42}. This estimate and the inequality (3.17) in \cite{ENV} imply that
$$
t\eta(\{|\widetilde{\mathfrak{R}}_{\eta,\e}^K f(x)|>t\}\le t\eta(\{|\widetilde{R}_{\nu,\e}^K(x)|>\tfrac t2\})\le2C\|\nu\|\le2C\|f\|_{L^1(\eta)},
$$
and we obtain (3.16) in \cite{ENV}. Now we complete the proof of Lemma 3.4 exactly as in \cite{ENV}. 
Theorem \ref{th24} follows directly from this Lemma 3.4 and Theorem 10.1 in \cite{NTV98}.
\end{proof}

\begin{proof}[Proof of Theorem \ref{th25}]
(i) The estimate \eqref{f211} is a corollary of Theorems \ref{th31} and \ref{th24}, and its proof is a repetition of the arguments in \cite{ENV}. Without loss of generality we assume that
$$
\mathbf S:=\sup_{x\in\supp\mu}\int_0^\infty\bigg[\frac{\mu(\BB(x,r))}{r^s}\bigg]^2\,\frac{dr}{r}<\infty.
$$
Otherwise (\ref{f32}) becomes trivial. We consider the measure
$$
\eta:=(2s\mathbf S)^{-1/2}\mu.
$$
Since for every $x\in\supp\mu$ and $r>0$,
$$
\mathbf S\ge\int_0^\infty\bigg[\frac{\mu(\BB(x,t))}{t^{s}}\bigg]^2\,\frac{dt}{t}\ge
[\mu(\BB(x,r))]^2\int_r^\infty\frac{dt}{t^{2s+1}}=\frac{[\mu(\BB(x,r))]^2}{2sr^{2s}},
$$
we see that $\eta\in\Sigma_s$. From \eqref{f31} we deduce
$$
\|\mathfrak{R}_{\eta,\e}^K\chi_Q\|^2_{L^2(\eta|Q)}=
(2s\mathbf S)^{-3/2}\|\mathfrak{R}_{\mu,\e}^K\chi_Q\|^2_{L^2(\mu|Q)}\le
C'\mathbf S^{-1/2}\mu(Q)=C''\eta(Q),
$$
where $C''=C''(A,s,\d)$. Thus, we are under the conditions of Theorem  \ref{th24}. By \eqref{f210},
$$
\|\mathfrak{R}_{\mu,\e}^K\|^2_{L^2(\mu)\to L^2(\mu)}=
(2s\mathbf S)\|\mathfrak{R}_{\eta,\e}^K\|^2_{L^2(\eta)\to L^2(\eta)}\le
C\mathbf S,\quad \e>0.
$$
The desired estimate \eqref{f211} follows immediately from this inequality.

(ii) The second part of Theorem \ref{th25} is a direct consequence of \eqref{f26}.
\end{proof}

\section{Relations between capacities}%Section 5

\begin{proof}[Proof of Theorem \ref{th26}] 
1. We start with the inequality
\begin{equation}\label{f51}
\g_{K,op}(E)\le C\g_{K,\ast}(E).
\end{equation}
Our arguments are similar to those in the proofs of Theorems 5.3 and 5.16 \cite[p.~29--30, 47]{V}. Let $\mu$ be a measure participating in \eqref{f213}. By \cite[Theorem 2.1]{NTV2002} or \cite[Theorem 5.13]{V} there exist absolute constants $\a>0$, $D>0$, and a function $h$ such that
$$
0\le h\le1,\quad\int_E h\,d\mu\ge\a\mu(E),\quad\|\mathfrak{R}_{\mu}^Kh\|_{L^\infty(\mu)}\le D.
$$
In spite of the fact that Nazarov, Treil and Volberg formulate this result for the Euclidean space, their proof works in our case without any changes. For $f\in L^1(\mu)$ set
$$
\widetilde{M}f(x):=\sup_{r>0}\frac1{\mu(\BB(x,3r))}\int_{\BB(x,r)}|f|\,d\mu.
$$
The Cotlar type inequality (see \cite[Theorem 7.1]{NTV98})
$$
\mathfrak{R}_{\mu,\ast}^Kh(x)\le C'\widetilde{M}(\mathfrak{R}_{\mu}^Kh(x))
+C''[\widetilde{M}(|h|^2)(x)]^{1/2},\quad x\in\supp\mu,
$$
yields the estimate
\begin{equation}\label{f52}
\mathfrak{R}_{\mu,\ast}^Kh(x)\le C,\quad x\in\supp\mu.
\end{equation}
Here the constants $C',C'',C$ depend only on the CZ parameters $A,s,\d$ of $K$.

The next step is to prove that
\begin{equation}\label{f53}
\mathfrak{R}_{\mu,\ast}^Kh(x)\le C\ \text{ for any } x\in\XX
\end{equation}
with another $C=C(A,s,\d)$. Choose $\e>0$ and $x\in\XX\setminus\supp\mu$. If $\supp\mu\subset\BB(x,\e)$ then $\mathfrak{R}_{\mu,\e}^Kh(x)\equiv0$. Otherwise we set
$$
r=\inf\{\dist(x,y):y\in\supp\mu\setminus\BB(x,\e)\}=\dist(x,z)>0,\quad z\in\supp\mu.
$$
We have 
\begin{multline*}
|\mathfrak{R}_{\mu,\e}^Kh(x)|=|\mathfrak{R}_{\mu,r}^Kh(x)|\\
\le\bigg|\int_{\BB(z,2r)\setminus\BB(x,r)}K(x,y)h(y)\,d\mu(y)\bigg|+
\int_{\XX\setminus\BB(z,2r)}|K(x,y)-K(z,y)|h(y)\,d\mu(y)+
|\mathfrak{R}_{\mu,2r}^Kh(z)|.
\end{multline*}
The first term is bounded since $|K(x,y)|\le Ar^{-s}$ as $r\le\dist(x,y)$, and $\mu\in\Sigma_s$. The last one is bounded by \eqref{f52} since $z\in\supp\mu$. Finally,
\begin{multline*}
\int_{\XX\setminus\BB(z,2r)}|K(x,y)-K(z,y)|h(y)\,d\mu(y)\overset{\eqref{f28}}\le
A\int_{\XX\setminus\BB(z,2r)}\frac{\dist(x,z)^\d}{\dist(z,y)^{s+\d}}\,d\mu(y)\\
=A\int_{2r}^\infty\frac{r^\d}{t^{s+\d}}\,d\mu(\BB(z,t))\le
A(s+\d)\int_{2r}^\infty\frac{r^\d t^s\,dt}{t^{s+\d+1}}=C(A,s,\d),
\end{multline*}
and we obtain \eqref{f53}. Define the measure $\s$ by the equality $d\s(x)=C_1^{-1}h(x)\,d\mu(x)$, where $C_1=\max(C,1)$ ($C$ is the constant in \eqref{f53}). Then $\s\in\Sigma_s$, and by \eqref{f53},
$\mathfrak{R}_{\s,\ast}^K{\bf1}(x)\le1$, $x\in\XX$. Thus, $\s$ participates in \eqref{f213}. Hence,
$$
\g_{K,\ast}(E)\ge\|\s\|=C_1^{-1}\int h\,d\mu\ge\a_1\mu(E),\quad \a_1=\a_1(A,s,\d)>0,
$$
and we get \eqref{f51} with $C=\a_1^{-1}$. Note that we did not use here that $\XX$ is geometrically doubling.

2. Now we prove the inverse inequality
\begin{equation}\label{f54}
\g_{K,\ast}(E)\le C\g_{K,op}(E).
\end{equation}
Choose a measure $\mu$ participating in \eqref{f214}, and fix $\e>0$. Clearly,
$\mathfrak{R}_{\mu,\e}^K{\bf1}(x)\le1$, but $\mathfrak{R}_{\mu,\e}^K$ is not an operator with a CZ kernel. Again as in the proof of Theorem \ref{th24}, we consider the same function $\phi(t)$ and the operator $\widetilde{\mathfrak{R}}_{\mu,\e}^K$ with the CZ kernel $K_\e(x,y)$. For $\mu\in\Sigma_s$ and $f\in L^1(\mu)$ we have
\begin{equation}\label{f55}
|\mathfrak{R}_{\mu,\e}^Kf(x)-\widetilde{\mathfrak{R}}_{\mu,\e}^Kf(x)|\le C\widetilde{M}f(x),
\quad x\in\XX,
\end{equation}
where $C$ depends on the CZ parameters of $K$. In particular, for $f(x)={\bf1}$, \eqref{f55} implies the estimate
\begin{equation}\label{f56}
|\widetilde{\mathfrak{R}}_{\mu,\e}^K{\bf1}(x)|\le C,\quad x\in\XX.
\end{equation}
Now we apply the nonhomogeneous $Tb$ theorem \cite[Theorem 2.10]{HM09} in the particular case $b_1=b_2={\bf1}$. According to this theorem,
\begin{equation}\label{f57}
\|\widetilde{\mathfrak{R}}_{\mu,\e}^K\|_{L^2(\mu)\to L^2(\mu)}\le
C(\|\widetilde{\mathfrak{R}}_{\mu,\e}^K\|_{BMO_\kappa^2(\mu)}+P+1),
\end{equation}
where $C$ depends on the doubling constant $N$ and on the CZ parameters of $K_\e$ (these parameters are independent of $\e$!). Furthermore, $P$ is the smallest (or "almost smallest") constant such that $|\langle\widetilde{\mathfrak{R}}_{\mu,\e}^K\chi_Q,\chi_Q \rangle|\le
P\mu(\l Q)$ for all balls $Q$ and for some fixed constant $\l>1$. Here $\langle f,g\rangle=\int fg\,d\mu$. Since the kernel $K_\e$ is antisymmetric, $\langle\widetilde{\mathfrak{R}}_{\mu,\e}^K\chi_Q,\chi_Q \rangle=0$ for any measurable set $Q$. Hence, $P=0$. Moreover, by \eqref{f56} the BMO-norm in \eqref{f57} is bounded by a constant depending only on $A,s,\d$. Thus,
$$
\|\widetilde{\mathfrak{R}}_{\mu,\e}^K\|_{L^2(\mu)\to L^2(\mu)}\le C,\quad C=C(N,A,s,\d).
$$
The maximal operator $\widetilde{M}$ is bounded on $L^2(\mu)$ -- see \cite[Lemma 3.1]{NTV98}. Hence, the operators $\mathfrak{R}_{\mu,\e}^K$ and $\widetilde{\mathfrak{R}}_{\mu,\e}^K$ are bounded simultaneously, and their norms differ at most by $C$. Thus,
$$
\|\mathfrak{R}_{\mu,\e}^K\|_{L^2(\mu)\to L^2(\mu)}\le C,\quad C=C(N,A,s,\d)
$$
(note that $C$ is independent of $\e$). We conclude that $C^{-1}\mu$ participates in \eqref{f213}. So, we have \eqref{f54}, and Theorem \ref{th26} is proved.
\end{proof}

\begin{lemma}\label{le51} 
Let $\mu\in\Sigma_\f$ with $\f$ satisfying \eqref{f216}. Then
\begin{equation}\label{f58}
\mathfrak{R}_{\mu,\ast}^\f{\bf1}(x)\le \|\mathfrak{R}_{\mu}^\f{\bf1}\|_\infty+C,\quad x\in\R^d,\quad C=C(s,\L).
\end{equation}
\end{lemma}
\begin{proof}
The idea of proof is not new -- see \cite[Lemma 2]{Vi} or \cite[p.~47]{V}. Note that \eqref{f216} implies the existence of $\mathfrak{R}_{\mu}^\f{\bf1}(x)$ almost everywhere in $\R^d$ with respect to Lebesgue measure. We may assume that 
$\|\mathfrak{R}_{\mu}^\f{\bf1}(x)\|_\infty<\infty$ (otherwise \eqref{f58} is trivial). Fix $\e>0$ and $x\in\R^d$. Let $\LL^d$ be Lebesque measure in $\R^d$, and let $\e_1>0$ be such that $\psi(\e_1)=\frac12\psi(\e)$, where $\psi$ is the function defined in Lemma \ref{le32}. Consider the mean value integral
$$\aligned
\frac1{\a_d\e_1^d}\int_{B(x,\e_1)}\int_{B(x,\e)}\frac{d\mu(y)}{\f(|y-z|)}\,d\LL^d(z)&<
\frac1{\a_d\e_1^d}\int_{B(x,\e)}\int_0^{\e_1}\frac{\a_dt^{d-1}\,dt}{\f(t)}\,d\mu(y)\\
&<\frac{\L}{\e_1^d}\,\frac{\e_1^d}{\f(\e_1)}\,\mu(B(x,\e))<C_0,
\endaligned$$
where $\a_d$ is the Lebesque measure of the unit ball in $\R^d$. Hence, there is a point 
$z\in B(x,\e_1)$ such that $|\mathfrak{R}_{\mu}^\f{\bf1}(z)|\le 
\|\mathfrak{R}_{\mu}^\f{\bf1}\|_\infty$, and
$$
\int_{B(x,\e)}\frac{d\mu(y)}{\f(|y-z|)}\le C_0.
$$
We have
$$\aligned
|\mathfrak{R}_{\mu,\e}^\f{\bf1}(x)|&\le|\mathfrak{R}_{\mu}^\f{\bf1}(z)|+
|\mathfrak{R}_{\mu,\e}^\f{\bf1}(x)-\mathfrak{R}_{\mu}^\f{\bf1}(z)|\\
&\le\|\mathfrak{R}_{\mu}^\f{\bf1}\|_\infty+\int_{\R^d\setminus B(x,\e)}
|K_\f(x,y)-K_\f(z,y)|\,d\mu(y)+\int_{B(x,\e)}|K_\f(z,y)|\,d\mu(y).
\endaligned$$
Since $\psi(|x-z|)\le\psi(\e_1)=\frac12\psi(\e)\le\frac12\psi(|x-y|)$, $y\in\R^d\setminus B(x,\e)$, we may apply the property \eqref{f28} of $K_\f$ (see the part (ii) of Lemma \ref{le32}). Using \eqref{f38} and integrating by parts, we get
\begin{multline}\notag
\int_{\R^d\setminus B(x,\e)}|K_\f(x,y)-K_\f(z,y)|\,d\mu(y)\le
C(s)\int_{\R^d\setminus B(x,\e)}\frac{\f(|x-z|)^{1/s}}{\f(|x-y|)^{1+1/s}}\,d\mu(y)\\
<C(s)\int_\e^\infty\frac{\f(\e)^{1/s}}{\f(t)^{1+1/s}}\,d\mu(B(x,t))
\le C'(s)\f(\e)^{1/s}\int_\e^\infty\frac{d\f(t)}{\f(t)^{1+1/s}}=sC'(s).
\end{multline}
Thus, we have \eqref{f58}, and Lemma \ref{le51} is proved.
\end{proof}

\begin{proof}[Proof of Theorem \ref{th27}] We start with the first inequality in \eqref{f217}. Let $\mu$ be a measure participating in \eqref{f214}. Then $|\mathfrak{R}_{\mu}^\f{\bf1}(x)|\le1$ $\LL^d$-a.~e. in $\R^d$. Moreover, $\Sigma_s\subset\Sigma_\f$ (see the part (iii) of Lemma \ref{le32}). Thus, $\mu$ participates in the definition of $\g_{\f,+}$ (see Section~1), and we have the desired inequality.

The second inequality in \eqref{f217} is a direct consequence of Lemma \ref{le51}.
\end{proof}

\section{Proof of Theorem \ref{th21} and Corollary \ref{cor22}}%Section 6

We prove a stronger assertion than the first part of Theorem \ref{th21}.

\begin{theorem}\label{th61}
 Let $\XX$ be a compact Hausdorff geometrically doubling metric space, and let $K$ be an antisymmetric CZ kernel. Then for every bounded closed set $E\subset\XX$,
\begin{equation}\label{f61}
\g_{K,\ast}(E)\ge c\,\sup\|\mu\|^{3/2}\bigg[\int_{\XX}W_s^\mu(x)\,d\mu(x)\bigg]^{-1/2},
\end{equation}
where the supremum is taken over all positive Radon measures supported by $E$, and $c$ depends only on the parameters of $K$ and on $N$.
\end{theorem}
\begin{proof}
In fact, the proof is a minor and obvious modification of the arguments in the proof of Theorem 2.7 in \cite{ENV}. Namely, one should replace $\g_{s,+}$ with $\g_{K,\ast}$, and use \eqref{f215} instead of (10.2) in \cite{ENV}. We omit details.
\end{proof}

\begin{proof}[Proof of Theorem \ref{th21}]
Since $\g_{\f,+}(E)\ge\g_{K_\f,\ast}(E)$, the first part is a corollary of Theorem~\ref{f61}.

To prove the second part, we rewrite \eqref{f213} in the form
\begin{equation}\label{f62}
\g_{\f,op}(E)=\sup\{\k\|\mu\|:\k>0,\ \k\mu\in\Sigma_s,\ \supp\mu\subset E,\
\pmb|\mathfrak{R}_{\k\mu}^\f\pmb|\le1\}.
\end{equation}
For any $\k$ and $\mu$ participating in the right hand side of \eqref{f62} we have
$$
1\ge\pmb|\mathfrak{R}_{\k\mu}^\f\pmb|^2=\k^2\pmb|\mathfrak{R}_{\mu}^\f\pmb|^2
\overset{\eqref{f212}}\ge\frac{c\k^2}{\|\mu\|}\int_{\R^d}W_\f^\mu(x)\,d\mu(x).
$$
Hence,
$$
\k\le C\|\mu\|^{1/2}\bigg[\int_{\R^d}W_{\f}^\mu(x)\,d\mu(x)\bigg]^{-1/2},\quad C=C(s).
$$
From \eqref{f62} we obtain the estimate
$$
\g_{\f,op}(E)=\sup\k\|\mu\|\le C\sup\|\mu\|^{3/2} \bigg[\int_{\R^d}W_\f^\mu(x)\,d\mu(x)\bigg]^{-1/2}.
$$
To deduce \eqref{f23} from this estimate and from \eqref{f215},  \eqref{f217}, we have to check only that the conditions on $\f$ imply \eqref{f216}. Indeed,
$$\aligned
\int_0^r\frac{t^{d-1}\,dt}{\f(t)}&=
\sum_{i=0}^\infty\int_{2^{-i-1}r}^{2^{-i}r}\frac{t^{d-1}\,dt}{\f(t)}
\le\sum_{i=0}^\infty\frac{2^{-di}r^d}{\f(2^{-(i+1)}r)}\\
&\overset{\eqref{f21}}\le\sum_{i=0}^\infty\frac{2^{-di}r^d}{2^{-(i+1)s}\f(r)}
=2^s\frac{r^d}{\f(r)}\sum_{i=0}^\infty2^{-(d-s)i}=C(d,s)\frac{r^d}{\f(r)}.
\endaligned$$
Theorem \ref{th21} is proved.
\end{proof}

\begin{proof}[Proof of Corollary \ref{cor22}]
The proof is based on the following inequality of Wolff, see \cite[p.~109, Theorem~4.5.2]{AH}:
for $1<p<\infty,\ \a p\le d$,
\begin{equation}\label{f63}
\|G_\a\ast\mu\|_{p'}^{p'}\approx\int_{\R^d}W_{\a,p}^\mu(x)\,d\mu(x),\quad
W_{\a,p}^\mu(x)=\int_0^1\biggl(\frac{\mu(B(x,t))}{t^{d-\a p}}\biggr)^{p'-1}\frac{dt}{t}.
\end{equation}
Take $\a=\frac23d,\ p=\frac32$. Then $p'=3,\ d-\a p=0$, and
\begin{equation}\label{f64}
W_{\frac23d,\frac32}^\mu(x)=\int_0^1\mu(B(x,t))^2\frac{dt}{t}.
\end{equation}
Clearly,
$$\aligned
\int_0^{e^{-3/2}}\biggl(\frac{\mu(B(x,t))}{\f(t)}\biggr)^{2}\frac{d\f(t)}{\f(t)}&=
\int_0^{e^{-3/2}}\mu(B(x,t))^2\biggl(\log\frac{1}{t}\biggr)^{3/2}
d\biggl(\log\frac{1}{t}\biggr)^{-1/2}\\
&=\frac12\int_0^{e^{-3/2}}\mu(B(x,t))^2\frac{dt}{t}.
\endaligned$$
Since diam$(E)\le1$, we may cover $E$ by at most $N=N(d)$ balls of radius $\frac12e^{-5/2}$. Let $B'$ be the ball with maximal measure. Then
\begin{equation}\label{f65}\begin{split}
\int_{\R^d}\int_0^{e^{-3/2}}\mu(B(x,t))^2\frac{dt}{t}\,d\mu(x)&\ge
\int_{B'}\int_0^{e^{-3/2}}\mu(B(x,t))^2\frac{dt}{t}\,d\mu(x)\\
&>\mu(B')^3>c(d)\|\mu\|^3.
\end{split}\end{equation}
Using this inequality we obtain the following estimates:
\begin{multline}\notag
\int_{\R^d}W_{\f}^\mu(x)\,d\mu(x)
\le\int_{\R^d}\biggl[\int_0^{e^{-3/2}}\mu(B(x,t))^2\frac{dt}{t}
+\int_{e^{-3/2}}^\infty\|\mu\|^2\frac{d\f(t)}{\f(t)^3}\biggr]\,d\mu(x)\\
=\int_{\R^d}\int_0^{e^{-3/2}}\mu(B(x,t))^2\frac{dt}{t}\,d\mu(x)+\frac34\|\mu\|^3
\overset{\eqref{f64},\eqref{f65}}<C(d)\int_{\R^d}W_{\frac23d,\frac32}^\mu(x)\,d\mu(x).
\end{multline}
On the other hand,
$$\aligned
\int_{\R^d}W_{\frac23d,\frac32}^\mu(x)\,d\mu(x)
&\le\int_{\R^d}\biggl[\int_0^{e^{-3/2}}\mu(B(x,t))^2\frac{dt}{t}
+\|\mu\|^2\int_{e^{-3/2}}^1\frac{dt}{t}\biggr]\,d\mu(x)\\
&\le C'(d)\int_{\R^d}W_{\f}^\mu(x)\,d\mu(x).
\endaligned$$
Thus,
 \begin{equation}\label{f66}
\|G_{\frac23d}\ast\mu\|_{3}^{3}\overset{\eqref{f63}}\approx
\int_{\R^d}W_{\frac23d,\frac32}^\mu(x)\,d\mu(x)\approx\int_{\R^d}W_{\f}^\mu(x)\,d\mu(x).
\end{equation}
Choose $\mu$ for which
$$
C_{\frac23d,\frac32}(E)\le2\|\mu\|^{3/2}\|G_{\frac23d}\ast\mu\|_{3}^{-3/2}.
$$
By \eqref{f22} and \eqref{f66} we have
$$
\g_{\f,+}(E)\ge
c\,\|\mu\|^{3/2}\|G_{\frac23d}\ast\mu\|_{3}^{-3/2}\ge c'C_{\frac23d,\frac32}(E).
$$

To prove the inverse inequality we note that $\f(t)$ satisfies the conditions (ii) of Theorem~\ref{th21}. By \eqref{f23} there is a measure $\mu\in M_+(E)$ for which
$$
\g_{\f,+}(E)\le2C\|\mu\|^{3/2} \bigg[\int_{\R^d}W_\f^\mu(x)\,d\mu(x)\bigg]^{-1/2},
$$
where $C$ is the constant in \eqref{f23}. Then
$$
C_{\frac23d,\frac32}(E)\ge\|\mu\|^{3/2}\|G_{\frac23d}\ast\mu\|_{3}^{-3/2}
\overset{\eqref{f66}}\approx\|\mu\|^{3/2} \bigg[\int_{\R^d}W_\f^\mu(x)\,d\mu(x)\bigg]^{-1/2}
\ge c\,\g_{\f,+}(E),
$$
and the proof is completed.
\end{proof}
\vspace{.2cm}
\noindent{\bf Acknowledgement.}
The authors are grateful to Professor F.~Nazarov for valuable suggestions and discussions.

%\begin{theorem}\label{th} \end{theorem} \begin{equation}\label{f}\end{equation} \overset{\eqref{f}}>\begin{equation}\begin{split}\end{split}$$\aligned\endaligned$$ \eqref{f} \begin{multline}\notag\end{multline}\begin{multline}\end{multline}

\end{document}